\documentclass[11pt]{article}
\usepackage{amsmath,amsthm,amssymb,graphicx,wrapfig}

\newtheorem{theorem}{Theorem}[section]

\newtheorem{proposition}[theorem]{Proposition}

\newtheorem{lemma}[theorem]{Lemma}

\newtheorem{corollary}[theorem]{Corollary}
\newtheorem{example}[theorem]{Example}

\theoremstyle{remark}

\newcommand{\sign}{\mathop{\mathrm{sign}}\nolimits}

\newcommand{\R}{\mathbb{R}}
\numberwithin{equation}{section}
\title{The Rasmussen invariant 
of a homogeneous knot}
\date{}
\author{Tetsuya Abe \\ 
\footnotesize{Department of Mathematics, Osaka City University} \\[-3pt]
\footnotesize{Sugimoto, Sumiyoshi-ku Osaka 558-8585, Japan}\\
\footnotesize{Email: t-abe@sci.osaka-cu.ac.jp}}

\begin{document}

\maketitle
 
\begin{abstract}
A homogeneous knot is a generalization of 
alternating knots and positive knots. 
We determine the Rasmussen invariant of a homogeneous knot.
This is a new class of knots such that the Rasmussen invariant is explicitly
described in terms of its diagrams.
As a corollary, 
we obtain some characterizations of a positive knot.
In particular, 
we recover Baader's theorem which states 
that a knot is positive if and only if 
it is homogeneous and strongly quasipositive.
\end{abstract}

\section{Introduction}
In \cite{Rasmussen}, Rasmussen  introduced a smooth
concordance invariant of a knot $K$ by using the Khovanov-Lee theory (see 
\cite{Khovanov-original} and \cite{Lee}),
now called the Rasmussen invariant $s(K)$.
This gives a lower bound for 
the four ball genus $g_{*}(K)$ of a knot $K$ as follows. 
\begin{equation} \label{eq:lower-bound}
 |s(K)| \le 2g_{*}(K).
\end{equation}
This lower bound is very powerful and it enables us to
 give a combinatorial proof of the Milnor conjecture
on the unknotting number of a torus knot.
Our motivation for studying the Rasmussen invariant is 
to describe $s(K)$ in terms of a given diagram of a knot $K$
to better understand $g_{*}(K)$.
From this point of view, some estimations of the Rasmussen invariant
are known (Plamenevskaya \cite{Plamenevskaya}, Shumakovitch \cite{Shumakovitch} and Kawamura \cite{Kawamura}.
See also Stoimenow \cite{Stoimenow}).



Let $O_{+}(D)$\footnote{In \cite{Kawamura2} and \cite{Lobb}, 
it was denoted by $l_{0}(D)$ and $\#$components$(T^{+}(D))$ respectively.} 
and $O_{-}(D)$ be the numbers of connected components of the diagrams 
which is obtained from $D$
by smoothing all negative and positive crossings of $D$, respectively.
Recently, Kawamura \cite{Kawamura2} and Lobb \cite{Lobb} independently 
obtained a more sharper estimation for the Rasmussen invariant
as follows.
	
\begin{theorem} [\cite{Kawamura2} and \cite{Lobb}] \label{theorem:Kawamura-Lobb}
Let $D$ be a  diagram of a knot $K$. Then
\[ w(D)-O(D)+2O_{+}(D)-1 \le s(K),\]
where $\omega (D)$ denotes the writhe of $D$ 
(i.e.~the number of positive crossings of $D$ minus the number of negative crossings of $D$) 
and $O(D)$ denotes the number of the Seifert circles of $D$.
\end{theorem}
Let $\Delta(D)=O(D)+1-O_{+}(D)-O_{-}(D)$
(a graph theoretical interpretation of $\Delta(D)$ 
due to Lobb is given 
in Section \ref{section:A characterization of a homogeneous knot}).	
In addition to Theorem \ref{theorem:Kawamura-Lobb}, 
Lobb \cite{Lobb} showed that if $\Delta(D)=0$,
then $s(K)= w(D)-O(D)+2O_{+}(D)-1$.

Our motivation for this paper is to study 
which diagrams $D$ satisfy the condition $\Delta(D)=0$.
Lobb \cite{Lobb} showed that 
if $D$ is positive, negative, alternating, 
or a certain braid diagram,
then $\Delta(D)=0$.
Note that these diagrams are all homogeneous
(the definition is given in Section 2). 
In this paper, 
we show that if $D$ is a homogeneous diagram of a knot,
then $\Delta(D)=0$ (the converse is also true. See Theorem \ref{theorem:homog})
and our main result is 
to determine the Rasmussen invariant of a  homogeneous knot.
This is a new class of knots such that the Rasmussen invariant is
 explicitly described in terms of its diagrams.
\begin{theorem} \label{theorem:Rasmussen-homog}
Let $D$ be a  homogeneous diagram of a knot $K$.
Then 
\[s(K)= w(D)-O(D)+2O_{+}(D)-1.\]
\end{theorem}

Ozsv\'{a}th and Szab\'{o} \cite{Oz2} and 
Rasmussen \cite{Rasmussen2} independently introduced another smooth
concordance invariant of a knot $K$ by using the Heegaard Floer homology theory,
now widely known as the tau invariant $\tau(K)$.
The Rasmussen invariant and tau invariant share some formal properties
and these are closely related 
to positivity of knots.
There are many notions of positivity (e.g.~braid positive, positive,
strongly quasipositive and quasipositive).
We recall these notions of positivity in Section 3.
Let $D$ be a diagram of a knot $K$.
Then Kawamura \cite{Kawamura2} also proved
\[ w(D)-O(D)+2O_{+}(D)-1 \le 2\tau(K).\]
Note that, if $\Delta(D)=0$,
$2 \tau(K)= w(D)-O(D)+2O_{+}(D)-1$.\footnote{
by using the fact that $-\tau(K)=\tau(\overline{K})$ for any knot $K$ 
\cite{Oz2}, 
where $\overline{K}$ denotes the mirror image of $K$.}
Therefore the corresponding result to Theorem
\ref{theorem:Rasmussen-homog}
holds for the tau invariant. 
In particular, we obtain $\tau(K)=s(K)/2 $ for a homogeneous knot $K$.

On the other hand, 
the Rasmussen invariant and tau invariant sometimes behave differently.
It has been conjectured that $\tau=s/2$, however, Hedden and Ording \cite{Hedden} proved that the Rasmussen invariant and tau invariant
are distinct (see also \cite{Liv2}).
It may be worth remarking that 
the Rasmussen invariant is sometimes stronger than the tau invariant
as an obstruction to a knot being smoothly slice 
(\cite{Hedden} and \cite{Liv2}, see also \cite{Freedman}). 
This is the reason why we are more interested in the Rasmussen invariant
rather than the tau invariant.

One can easily see that a braid positive knot is strongly quasipositive,
however, it is not obvious whether a positive knot is strongly quasipositive.
Nakamura \cite{Nakamura} and Rudolph \cite{Rudolph} independently
proved that a positive knot is strongly quasipositive.
Not all strongly quasipositive knots are positive. 
For instance, such examples are given by divide knots \cite{Rudolph2}.
Rudolph \cite{Rudolph} asked whether positive knots could be 
characterized  as strongly positive knots with some extra
geometric conditions. 
Several years later, Baader found that the extra condition is homogeneity.
To be precise, 
Baader \cite{Baader} proved that a knot is positive if and only
if it is homogeneous and strongly quasipositive.
As a corollary of Theorem \ref{theorem:Rasmussen-homog}, 
we obtain some characterizations of a positive knot.
\begin{theorem}  \label{theorem:positivity}
Let $K$ be a knot. Then $(1)$--$(4)$ are equivalent.\\
$(1)$  $K$ is  positive.\\
$(2)$  $K$ is  homogeneous and strongly quasipositive.\\
$(3)$  $K$ is  homogeneous, quasipositive and $g_{*}(K)=g(K)$.\\
$(4)$  $K$ is  homogeneous and $\tau(K)=s(K)/2=g_{*}(K)=g(K)$.
\end{theorem}
In particular, we recover Baader's theorem. 
Note that our proof is 4-dimensional in the sense that 
we use concordance invariants, whereas Baader \cite{Baader}
used the Homflypt polynomial.
As an immediate corollary of Theorem \ref{theorem:positivity},
we obtain the following.
\begin{corollary}  \label{corollary:positivity}
Let $K$ be a  homogeneous knot.
Then  the following are equivalent.\\
$(1)$ $K$ is  positive.\\ 
$(2)$ $K$ is  strongly quasipositive.\\ 
$(3)$ $K$ is  quasipositive and $g_{*}(K)=g(K)$.\\ 
$(4)$ $\tau(K)=s(K)/2=g_{*}(K)=g(K)$.
\end{corollary}

It may be interesting to compare Corollary \ref{corollary:positivity} 
and the following proposition by Hedden.
\begin{proposition} [\cite{Hedden2}]
Let $K$ be a fibered  knot.
Then  the following are equivalent.\\
$(1)$ $K$ is strongly quasipositive.\\
$(2)$ $K$ is quasipositive and $g_{*}(K)=g(K)$.\\
$(3)$ $\tau(K)=g_{*}(K)=g(K)$.
\end{proposition}
At a first glance, we wonder why similar results
hold for fibered  knots and  homogeneous knots.
However, it is not surprising since
homogeneous knots are related to fiberedness.
For instance,
a knot which admits a homogeneous braid diagram is fibered
(see Section \ref{section:Geometric aspect of a homogeneous knot}
or Proposition 1.4~in \cite{Lobb}).

This paper is constructed as follows.
In Section \ref{section:Geometric aspect of a homogeneous knot},
we observe a geometric aspect of a homogeneous knot.
In Section \ref{section:A characterization of a homogeneous knot},
we give a new characterization of a homogeneous diagram of a knot
and determine the Rasmussen invariant of a homogeneous knot
(Theorem  \ref{theorem:Rasmussen-homog}).
In Section \ref{section:Positivity of knots}, 
we recall some notions of positivity for knots
and give  some characterizations of a positive knot
(Theorem \ref{theorem:positivity}).
In Section \ref{section:A new approach to estimate the Rasmussen invariant of a knot},
we propose a new approach to estimate the Rasmussen invariant of a knot.

\section*{Acknowledgments}
The author would like to thank the members
of Friday Seminar on Knot Theory in Osaka City University,
especially Masahide Iwakiri and In Dae Jong.
The author would like to thank Seiichi Kamada, Tomomi Kawamura 
and Kengo Kishimoto for
helpful comments on an earlier draft of this paper
and Mikami Hirasawa for explaining to him cut-and-paste arguments of Seifert surfaces in detail.
This work was supported by Grant-in-Aid for JSPS Fellows.

\section{Geometric aspect of a homogeneous knot}
\label{section:Geometric aspect of a homogeneous knot}
Cromwell \cite{cromwell} introduced the notion of homogeneity for
knots  to
generalize results on alternating knots.
The notion of homogeneity is also defined  for
signed graphs and diagrams.
For graph theoretical terminologies in this paper, 
we refer the reader the book of Cromwell \cite{cromwell-book}.\footnote{In this paper, 
we use the notation ``cycle" instead of ``circuit".}

A graph is \textit{signed} if
each edge of the graph is labeled $+$ or $-$.
A typical signed graph is the Seifert graph $G(D)$ associated to 
a knot diagram $D$:
for each  Seifert circle of $D$, we associate a vertex of $G(D)$ and
two vertices of $G(D)$ are connected by an edge if there is  a crossing of $D$
whose adjacent two Seifert circles are corresponding to the two vertices.
Each edge of $G(D)$ is labeled $+$ or $-$ depending
on the sign of its associated crossing of $D$.
For convenience, 
we say a $+$ or $-$ edge instead of an edge labeled $+$ or $-$.

A \textit{block} of a (signed) graph is a maximal subgraph of the graph 
with no cut-vertices.
A  signed graph is \textit{homogeneous} 
if each block has the same signs.
A diagram $D$ of a knot is \textit{homogeneous} 
if $G(D)$ is homogeneous.
Cromwell \cite{cromwell} showed that alternating diagrams and positive diagrams
are homogeneous.
 There are many homogeneous diagrams which are non-alternating
and non-positive.
\begin{example}
Let $D$ be the non-alternating and non-positive diagram 
as in Figure \ref{fig:two-bridge}.
Then $G(D)$ is homogeneous (see Figure \ref{fig:8_12-1}). Therefore $D$ is a homogeneous diagram which is non-alternating and non-positive.
Note that $D$ is not minimal crossing diagram
(this is not used later).
\begin{figure}
\begin{center}
\includegraphics[scale=0.35]{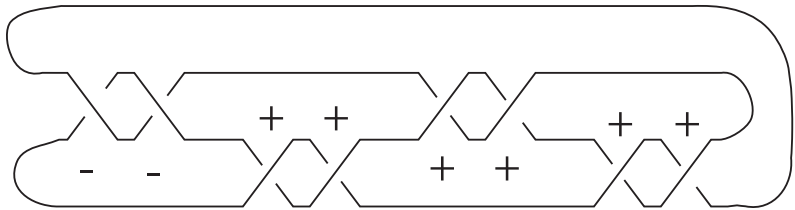}
\end{center}
\caption{a non-alternating and non-positive diagram}
\label{fig:two-bridge}
\end{figure}

\begin{figure}
\begin{center}
\includegraphics[scale=0.35]{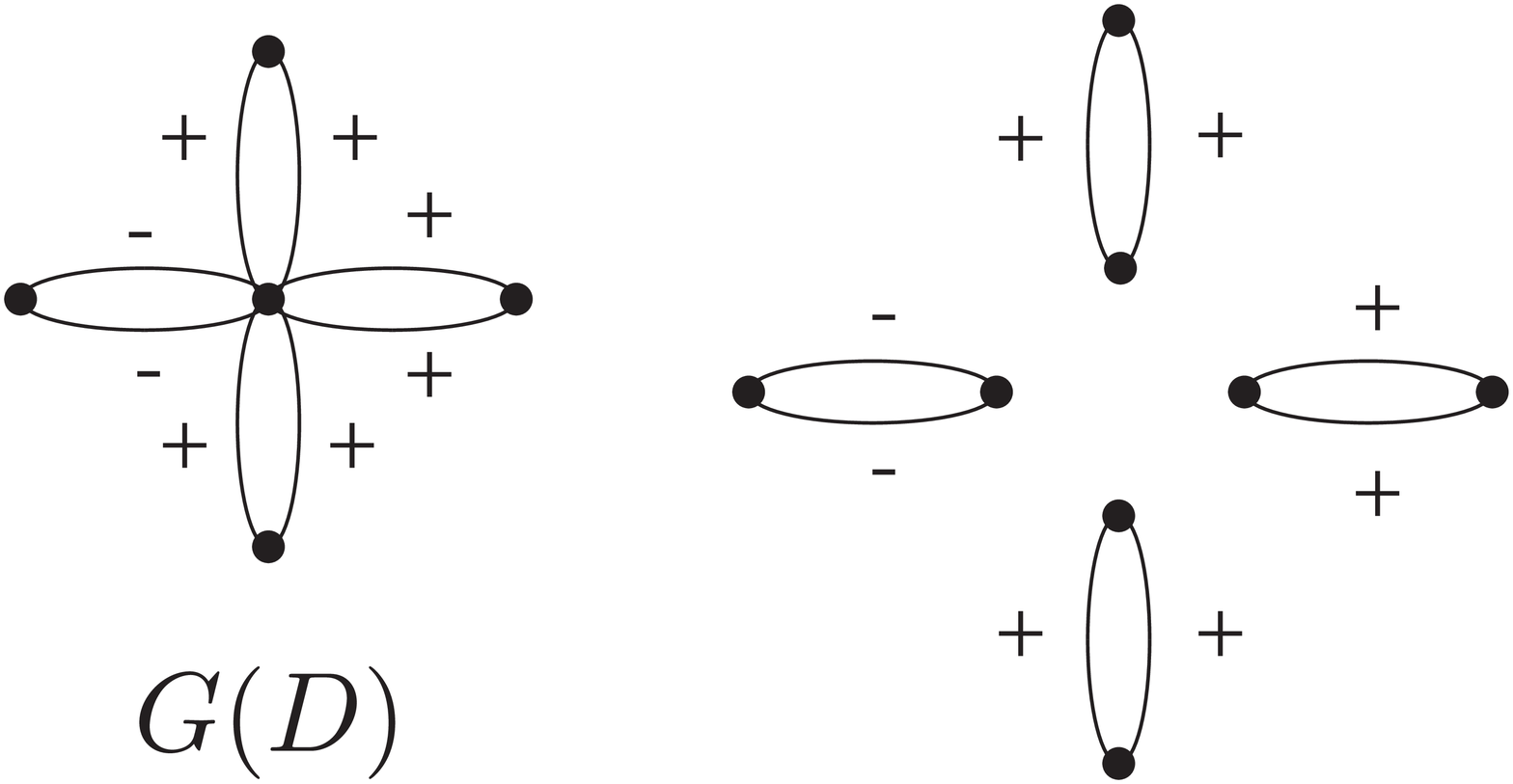}
\end{center}
\caption{the graph $G(D)$ is homogeneous}
\label{fig:8_12-1}
\end{figure}
\end{example}

Let $B_{n}$ be the braid group on $n$ strands 
with generators $\sigma_{1}, \sigma_{2},\cdots, \sigma_{n-1}$.
Stallings \cite{Stallings} 
introduced the notion of a homogeneous braid.
A braid $\beta=\sigma_{i_1}^{\epsilon_{1}} \sigma_{i_2}^{\epsilon_{2}}
\cdots \sigma_{i_k}^{\epsilon_{k}}, \epsilon_j=\pm1 \ (j=1, \cdots, k)$ is
\textit{homogeneous} if  
\begin{enumerate} 
\renewcommand{\labelenumi}{(\theenumi)}
\item  every $\sigma_{j}$ occurs at least once,
\item for each $j$, the exponents of all occurrences of $\sigma_{j}$
are the same.
\end{enumerate} 
For example, 
the braid $\sigma_1 \sigma^{-1}_2 \sigma_1 \sigma^{-1}_2$
is homogeneous, however, 
the braid $\sigma^{2}_1 \sigma_2 \sigma_1 \sigma^{-1}_2$ 
is not homogeneous. 
Stallings \cite{Stallings} proved that 
the closure of a homogeneous braid is fibered.
The following lemma is origin of the name ``homogeneous".
\begin{lemma} [\cite{cromwell}] 
Let $\beta$ be a braid whose closure is a knot.
Then $\beta$ is homogeneous if and only if 
the braid diagram of the closure of $\beta$ is homogeneous.
\end{lemma}
A knot $K$ is \textit{homogeneous} 
if $K$ has a homogeneous diagram.
The class of homogeneous knots includes alternating knots and positive knots.
There are homogeneous knots which are non-alternating and non-positive
and Cromwell \cite{cromwell} showed that the knot $9_{43}$ is the simplest one.
One of the distinguished properties of a homogeneous diagram is the following. 	
\begin{theorem} [\cite{cromwell}] \label{theorem:homogeneous}
Let $D$ be a homogeneous diagram of a knot $K$.
Then the genus of $K$ is realized by that of the Seifert surface obtained by applying Seifert's algorithm to $D$.
\end{theorem}
Cromwell proved the above theorem algebraically.
There is a geometric proof.
Here we give an outline of the proof,
which is suggested by M.~Hirasawa.

The Seifert circles of a diagram is divided into two types:
a Seifert circle is \textit{of type 1}
 if it does not contain any other Seifert circles in $\R ^{2}$,
otherwise it is \textit{of type 2}.
Let $D \subset \R^{2}$ be a knot diagram and $C$ 
 a type 2 Seifert circle of $D$.
Then $C$ separates $\R^{2}$ into two components $U$ and $V$ such that
$U \cup V =\R^{2}$ and $U \cap V = \partial U=\partial V=C$.
Let $D_1$ and $D_2$ be the diagrams formed form $D \cap U$ and 
$D \cap V$ by adding suitable arcs from $C$ respectively.
If both $(U-C) \cap D \neq \emptyset$ and $(V-C) \cap D \neq \emptyset$,
then $C$ \textit{decomposes} $D$ into
a $*$-product of $D_1$ and $D_2$, which is denoted by $D=D_1 * D_2$.
Then the Seifert surface obtained by applying Seifert's algorithm to $D$ is
 a Murasugi sum of Seifert surfaces obtained by applying Seifert's algorithm to 
$D_1$ and $D_2$ respectively 
(for the definition of a Murasugi sum, see \cite{Kaw} or \cite{Gabai}).
A diagram is \textit{special} if $D$ has no decomposing Seifert circles
of type 2.
A special positive (or negative) diagram is alternating.
Cromwell implicitly showed  the following (see Theorem 1 in \cite{cromwell}).
\begin{lemma} [\cite{cromwell}] \label{lemma:murasugi-sum}
Let $D$ be a homogeneous diagram of a knot $K$.
Then \\
(1) there are special diagrams $D_1, \cdots, D_n$ such that
$D=D_{1}*D_{2}*\cdots*D_{n}$,\\
(2) each special diagram $D_i \ (i=1, \cdots, n)$ is the connected sum of special alternating diagrams,\\
(3) each special alternating diagram corresponds to a block of $G(D)$.
\end{lemma}
Let $D$ be homogeneous diagram of a knot $K$.
Then, by Lemma \ref{lemma:murasugi-sum}, 
the Seifert surface $S$ obtained by applying Seifert's algorithm to $D$
 is Murasugi sums of the Seifert surfaces
obtained by applying Seifert's algorithm to the special alternating diagrams.
The following lemma is classical results of  Crowell and Murasugi.
\begin{lemma} \label{lemma:alternating}
Let $D$ be a alternating diagram of a knot $K$.
Then the genus of $K$ is realized by that of the Seifert surface obtained by applying Seifert's algorithm to $D$.
\end{lemma}
In \cite{Gabai2}, 
Gabai gave an elementary  proof of Lemma \ref{lemma:alternating}
by using cut-and-past arguments.
By Lemma \ref{lemma:alternating},
$S$ is Murasugi sums of minimal Seifert surfaces.
Let $R_1$ and $R_2$ be two minimal Seifert surfaces.
Then a Murasugi sum of $R_1$ and $R_2$ is 
a minimal Seifert surface due to Gabai \cite{Gabai}. 
Therefore we obtain a geometric proof of Theorem \ref{theorem:homogeneous}.  

\section{A characterization of a homogeneous diagram}
\label{section:A characterization of a homogeneous knot}
In this section, we give a new characterization of a homogeneous diagram 
(Theorem \ref{theorem:homog}).
In particular, we show that if $D$ is a homogeneous diagram of a knot,
then $\Delta(D)=0$. One can prove this by induction on the number of 
cut-vertices of $G(D)$, 
however, we prove this more graph-theoretically.
For the Seifert Graph $G(D)$ associated to a knot diagram $D$, 
we construct a graph (which is denoted by $G(D)_{\Delta}$ later)
such that the number of cycles of the graph is equal to $\Delta(D)$
and we prove if $D$ is homogeneous,
then the graph has no cycles.
By using Theorems \ref{theorem:Kawamura-Lobb} and \ref{theorem:homog}, 
we determine the Rasmussen invariant of a homogeneous knot
(Theorem  \ref{theorem:Rasmussen-homog}).

Let $G_{+}$ and $G_{-}$ be the graphs  
which are obtained from a signed graph $G$ by removing all $-$ and $+$ edges, 
respectively.
Here we note that, by definition, each vertex of $G$ belongs to 
exactly one connected component of $G_{+}$ and $G_{-}$ respectively. 
Let $G_{\Delta}$ be the graph 
whose vertices are the connected components of $G_{+}$ and $G_{-}$
and two vertices of $G_{\Delta}$ are connected by an edge
if a vertex of $G$  belong to the two connected components
(which correspond to the two vertices).  
We give two examples,
which provide us the idea of the proof of Lemma \ref{lemma:cycle}.

\begin{example}
Let $G(D)$ be the signed graph as in Figure \ref{fig:8_12-1}.
We label $1$ and $2$ the connected components of $G_{+}$ and
$3$, $4$, $5$ and $6$ the connected components of $G_{-}$.
Then $G_{\Delta}$ is the graph as in Figure \ref{fig:8_12-2}
and it is tree.
\begin{figure}
\begin{center}
\includegraphics[scale=0.35]{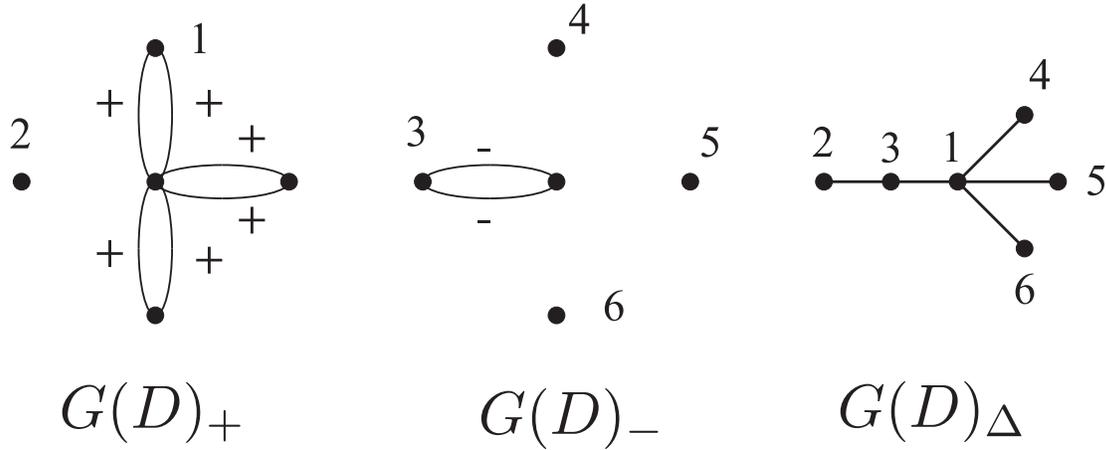}
\end{center}
\caption{the graph $G_{\Delta}$  is tree}
\label{fig:8_12-2}
\end{figure}
\end{example}
\begin{example}
Let $G$ be the signed graph as in Figure \ref{fig:cycle}.
Then $G$ has only one block and it is $G$ itself.
Since the block contains  $+$  and $-$ edges,
$G$ is not homogeneous. 
Note that $G$ has a cycle which contains  $+$  and $-$ edges
(in this case, the cycle is unique).
We label $1$ and $2$ the connected components of $G_{+}$ and
$3$, $4$, $5$ and $6$ the connected components of $G_{-}$.
Then $G_{\Delta}$ is the graph as in Figure \ref{fig:cycle}
and $G_{\Delta}$ has a cycle which is denoted by $(1,5)(5,2)(2,6)(6,1)$
(in this case, the cycle is also unique).

Conversely, let $\widetilde{e}_{1}, \widetilde{e}_{2},  \widetilde{e}_{3}$ 
and $\widetilde{e}_{4}$ be the edges of $G_{\Delta}$  
and $\widetilde{v}_{1}, \widetilde{v}_{2}, \widetilde{v}_{3}, \widetilde{v}_{4}$and $\widetilde{v}_{5}$  the vertices of $G_{\Delta}$
as in Figure \ref{fig:cycle2}.
Let $v_{i}$ $(i=1, \cdots, 4)$ be the vertex of $G$ which corresponds to 
$\widetilde{e}_{i}$  and let $v_{5}=v_{1}$.
Then $\widetilde{v}_{i+1}$ as a connected component of $G_{+}$ or $G_{-}$
contains $v_{i}$ and $v_{i+1}$
and there exists a simple path $l_{i}$  in $\widetilde{v}_{i+1}$ 
from $v_{i}$ to  $v_{i+1}$.
Therefore we obtain a cycle $l_{1} l_{2} l_{3} l_{4}$
from $v_{1}$ to  $v_{5} (=v_1)$.

\begin{figure}
\begin{center}
\includegraphics[scale=0.35]{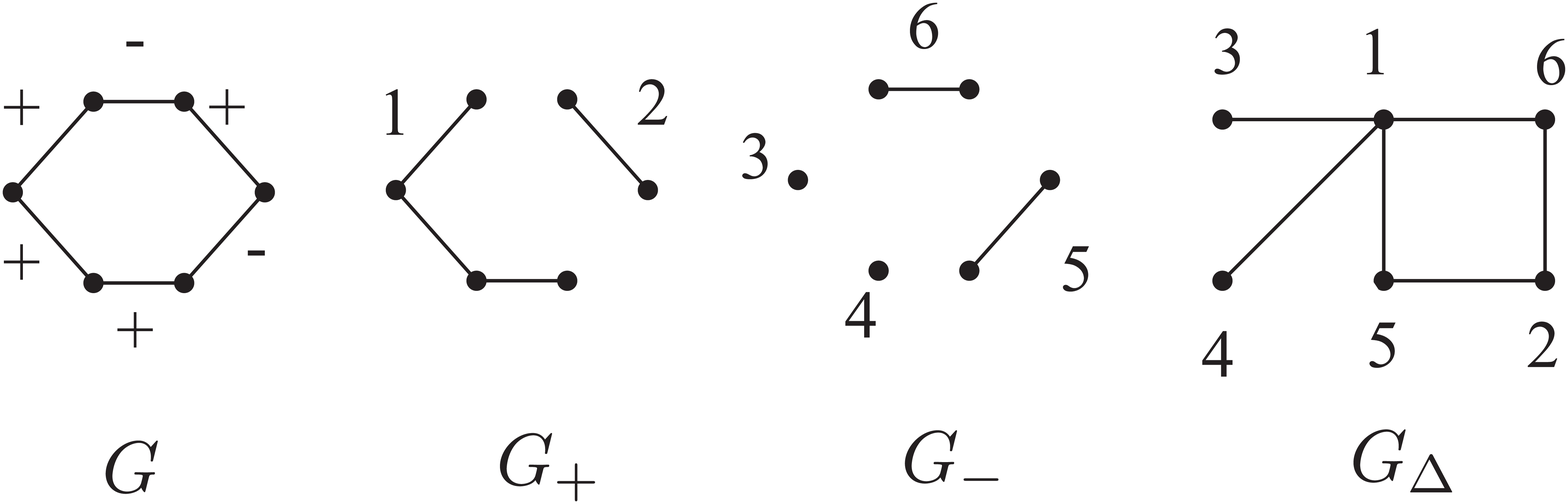}
\end{center}
\caption{the graph $G$ is not homogeneous}
\label{fig:cycle}
\end{figure}
\begin{figure}
\begin{center}
\includegraphics[scale=0.35]{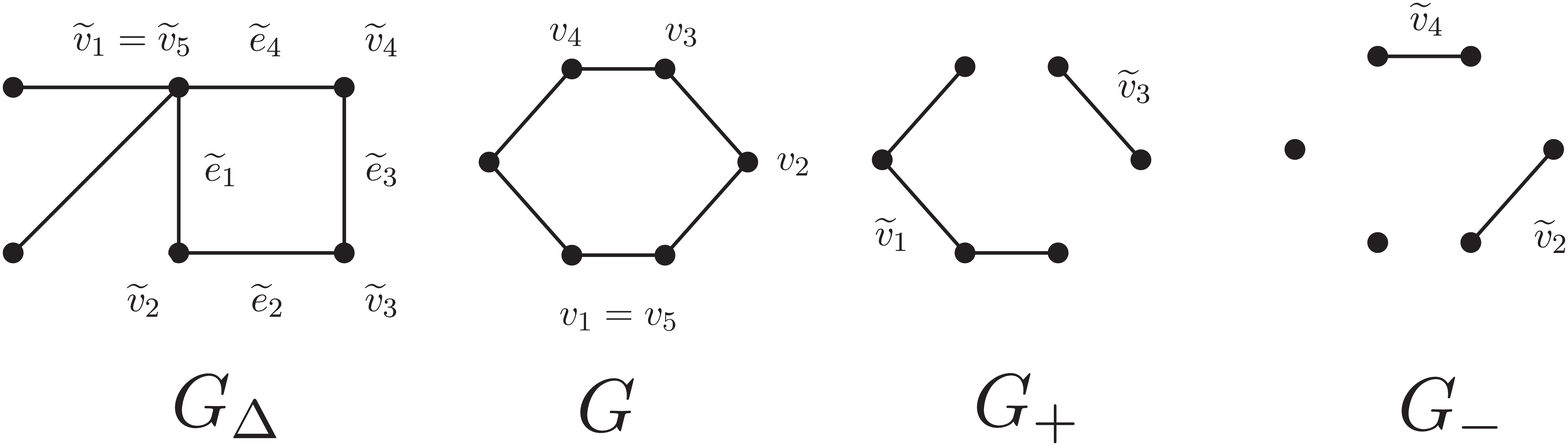}
\end{center}
\caption{the graph $G$ is not homogeneous}
\label{fig:cycle2}
\end{figure}
\end{example}
For a signed graph $G$, we denote by sign$(e)$ the sign of an edge $e$ of $G$.
We show the following lemma to prove Theorem \ref{theorem:homog}.
To prove $(2) \Longleftarrow (3)$ is essential.
\begin{lemma} \label{lemma:cycle}
Let $G$ be a signed graph. The following are equivalent.\\
(1) $G$ is not homogeneous.\\
(2) $G$ has a cycle which contains both $+$ and $-$ edges.\\
(3) $G_{\Delta}$ has a cycle. 
\end{lemma}

\begin{proof}
$(1) \Longrightarrow (2)$ 
Since $G$ is not homogeneous, by definition, there exists
a block which contains  $+$  and $-$ edges. 
Then there exist a vertex $v$ and edges $e_1$ and $e_2$ of the block
such that one of the endpoints of $e_1$ and $e_2$
is $v$ respectively and sign($e_1$) $\neq$ sign($e_2$). 
Now $v$ is not a cut-vertex
since $v$ is a vertex of the block.
There $G$ has a cycle which contains both $+$ and $-$ edges.\\
$(1) \Longleftarrow (2)$ 
Let $e_1  \cdots e_n$ be the cycle of $G$ which contains both $+$ and $-$ edges.
Then there exists a natural number $i$
such that $\sign(e_{i}) \neq \sign(e_{i+1})$.
Let $v$ be the vertex 
such that one of the endpoints of $e_i$ and $e_{i+1}$
is $v$ respectively. 
Since $v$ is not a cut-vertex,
edges $e_i$ and $e_{i+1}$ belong to the same block.
Therefore $G$ is not homogeneous.\\
$(2) \Longrightarrow (3)$ 
Let $(e_{1} \cdots e_{i_{1}}) 
\cdots
(e_{i_{k-1}+1} \cdots e_{i_{k}})
\cdots
(e_{i_{j-1}+1} \cdots e_{i_{j}})$ 
be the cycle which contains both $+$ and $-$ edges,
where $\sign (e_{i_{k-1}+1})=\cdots = \sign (e_{i_{k}})$ and
$\sign (e_{i_{k}}) \neq \sign(e_{i_{k}+1})$.
Then the path $e_{i_{k-1}+1} \cdots e_{i_{k}}$ is also one
in $G_{\sign(e_{i_{k}})}$, 
which contracts to a vertex $\widetilde{v}_{k}$ of $G_{\Delta}$
$(k=1, \cdots,j)$.
Let $\widetilde{e}_{k}$ be the edge of $G_{\Delta}$
whose endpoints are  
$\widetilde{v}_{k}$ and $\widetilde{v}_{k+1}$ $(k=1, \cdots, j-1)$,
which is corresponding to the vertex $v_{k}$  
such that one of the endpoints of $e_{i_{k}}$ and $e_{i_{k}+1}$ is $v_k$ respectively.
Let $\widetilde{e}_{j}$  be the edge of $G_{\Delta}$
whose endpoints are  
$\widetilde{v}_{j}$ and $\widetilde{v}_{1}$,
which is corresponding to the vertex $v_{j}$  
such that one of the endpoints of $e_{i_{j}}$ and $e_{1}$ is $v_j$ respectively.

Therefore $\widetilde{e}_{1} \cdots \widetilde{e}_{j}$
is a path from $\widetilde{v}_{1}$ to $\widetilde{v}_{1}$,
possibly not a cycle.
If the path is not a cycle, we choose a cycle 
(as a subsequence of edges of the path). 
Therefore $G_{\Delta}$ has a cycle.\\
$(2) \Longleftarrow (3)$
Let 
$\widetilde{e}_{1} \cdots \widetilde{e}_{n}$
be the cycle $G_{\Delta}$ $(i=1, \cdots, n)$ and denote
$\widetilde{e}_{i}=(\widetilde{v}_{i},\widetilde{v}_{i+1})$.
Then $\widetilde{v}_{n+1}=\widetilde{v}_{1}$.
Let $v_{i}$ be the vertex of $G$ which corresponds to 
$\widetilde{e}_{i}$ $(i=1, \cdots, n)$ and $v_{n+1}=v_{1}$.
Recall that a vertex  of $G_{\Delta}$ corresponds to
a connected component of $G_{+}$ or $G_{-}$.
Then $\widetilde{v}_{i+1}$ (as a connected component of $G_{+}$ or $G_{-}$)
contain $v_{i}$ and $v_{i+1}$ $(i=1, \cdots, n)$.
There exists a simple path $l_{i}$ from $v_{i}$ to  $v_{i+1}$.
There we obtain a path $l_{1} l_{2}\cdots l_{n}$
from $v_{1}$ to  $v_{n+1} (=v_1)$,
possibly not a cycle.
If the path is not a cycle, we choose a cycle (as a subsequence of edges of the path). By the construction, the cycle always 
contains both $+$ and $-$ edges.
\end{proof}

Note that $O_{+}(D)$ and $O_{-}(D)$ 
are equal to the numbers of connected components
of $G(D)_{+}$ and $G(D)_{-}$, respectively.
Therefore the number of vertices of $G(D)_{\Delta}$
is equal to $O_{+}(D) + O_{-}(D)$ and, by definition,
the number of edges of $G(D)_{\Delta}$
is equal to $O(D)$.
Lobb \cite{Lobb} showed that 
$\Delta(D)= b_{1}(G(D)_{\Delta})$ for any diagram $D$.
For the completeness, we recall the proof here. 
\begin{align*}
b_{1}(G(D)_{\Delta})&=b_{0}(G(D)_{\Delta})-\chi(G(D)_{\Delta})\\
                      &= 1 - (O_{+}(D) + O_{-}(D)-O(D))\\
                      &= \Delta(D),
\end{align*}
where $b_{i}$ denotes the $i$-th Betti number $(i=0,1)$ 
and $\chi$ denotes the Euler characteristic.
Then we obtain the following.
\begin{theorem} \label{theorem:homog}
A diagram $D$ of a knot is homogeneous 
if and only if $\Delta(D)=0$.
\end{theorem}

\begin{proof}
By the above argument, $\Delta(D)=0$ if and only if $G_{\Delta}$ is tree.
Therefore the proof immediately follows from Lemma \ref{lemma:cycle}.
\end{proof}

Now we prove Theorem \ref{theorem:Rasmussen-homog}.
\begin{proof}[Proof of Theorem \ref{theorem:Rasmussen-homog}]
By Theorem \ref{theorem:homog},
we obtain $\Delta(D)=0$.
As mentioned before, Lobb \cite{Lobb} showed that if $\Delta(D)=0$,
then $s(K)= w(D)-O(D)+2O_{+}(D)-1$. 
This completes the proof.
\end{proof}

\section{Positivity of knots}
\label{section:Positivity of knots}
In this section, we recall some notions of positivity
and give  some characterizations of a positive knot
(Theorem \ref{theorem:positivity}).
In particular, 
we recover Baader's theorem which states 
that a knot is positive if and only if 
it is homogeneous and strongly quasipositive.

\begin{figure}
\begin{center}
\includegraphics[scale=0.3]{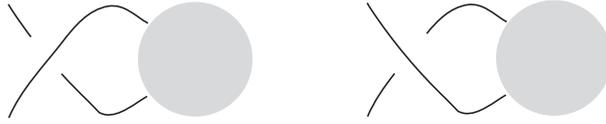}
\end{center}
\caption{nugatory crossings}
\label{fig:nugatory}
\end{figure}

Let $D$ be a diagram of a knot.
We denote by $D_{p}$ the diagram which is obtained from
$D$ by smoothing (along the orientation of $D$) at a crossing $p$.
A crossing of $D$ is \textit{nugatory} 
if there exists a curve $l$ such that 
the intersection of $D$ and $l$ is only a crossing of $D$
(see also Figure \ref{fig:nugatory}).
Then it is easy to see that the following lemma holds.
\begin{lemma} \label{lemma:nugatory}
Let $p$ be a crossing of $D$.
Then $p$ is nugatory 
if and only if 
the number of the connected components of $D_{p}$ is two.
\end{lemma}

Here we recall some notions of positivity for knots.
A knot is \textit{braid positive} 
if it is the closure of a braid of the form
$\beta =\prod_{k=1}^{m} \sigma_{i_{k}}$.  
A knot is \textit{positive} if it has  a diagram 
without negative crossings.
L.~Rudolph introduced the concept of a (strongly) quasipositive  knot
(see \cite{Rudolph2}).
Let 
\[ \sigma_{i,j} =(\sigma_{i}, \cdots, \sigma_{j-2})(\sigma_{j-1})
(\sigma_{i}, \cdots, \sigma_{j-2})^{-1}.
\]
A knot is \textit{strongly quasipositive} 
if it is the closure of a braid of the form
\[ \beta =\prod_{k=1}^{m} \sigma_{i_{k},j_{k}}.\]
A knot is \textit{quasipositive} 
if it is the closure of a braid of the form
\[ \beta =\prod_{k=1}^{m} \omega_{k} \sigma_{i_{k}} \omega_{k}^{-1},\]
where $\omega_{k}$ is a word in $B_{n}$.
The following are known.
\begin{enumerate} 
\renewcommand{\labelenumi}{(\theenumi)}
\item  Let $K$ be a torus knot.\footnote{A torus knot is strongly quasipositive.}
Then $\tau(K)=s(K)/2=g_{*}(K)=g(K)$, 
where $g(K)$ denotes the (Seifert) genus of $K$.
This is due to Rasmussen for $s$ \cite{Rasmussen}
and Ozsv\'{a}th and Szab\'{o} for $\tau$ \cite{Oz2}.
These equalities provide a proof of the Milnor conjecture.
\item Let $K$ be a strongly quasipositive knot. \label{page:1}
Then $\tau(K)=s(K)/2=g_{*}(K)=g(K)$.
This is due to Livingston \cite{Liv1}.
\item Let $K$ be a quasipositive knot. 
Then $\tau(K)=s(K)/2=g_{*}(K)$.
This is due to Plamenevskaya \cite{Plamenevskaya2}
 and Hedden  (with a detailed and constructive proof) \cite{Hedden2}
for $\tau$, and
Plamenevskaya \cite{Plamenevskaya} and Shumakovitch \cite{Shumakovitch} for $s$.\end{enumerate} 
By using Lemma \ref{lemma:nugatory}, we prove Theorem \ref{theorem:positivity}.
\begin{proof}[Proof of Theorem \ref{theorem:positivity}]
$(1) \Longrightarrow (2)$ A positive knot is strongly quasipositive
(\cite{Nakamura} and \cite{Rudolph}).\\
$(2) \Longrightarrow (3)$
A strongly quasipositive knot $K$ is a quasipositive knot with
$g_{*}(K)=g(K)$ \cite{Rudolph2}.\\
$(3) \Longrightarrow (4)$ 
Since $K$ is a quasipositive knot,
$\tau(K)=s(K)/2=g_{*}(K)$.
By the assumption, $g_{*}(K)=g(K).$
Therefore $\tau(K)=s(K)/2=g_{*}(K)=g(K)$.\\
$(4) \Longrightarrow (1)$
Let $D$ be a homogeneous diagram of $K$.
Then the genus of $K$ is realized by 
that of the surface constructed 
by applying Seifert's algorithm to $D$ (Theorem \ref{theorem:homogeneous}).
Therefore $2g(K)=1+c(D)-O(D)$,
where $c(D)$ denotes the number of crossings of $D$. 
By Theorem \ref{theorem:Rasmussen-homog}, we have
$s(K)= w(D)-O(D)+2O_{+}(D)-1$.
By assumption, $s(K)=2g(K)$.
This implies that $O_{+}(D)-1=n_{-}(D)$,
where $n_{-}(D)$ denotes the number of negative crossings of $D$.

If there exists a non-nugatory negative crossing $p$ of $D$,
then $D_{p}$ is connected by Lemma \ref{lemma:nugatory}.
Therefore $O_{+}(D) -1< n_{-}(D)$
(since, in general, 
the difference of the numbers of the connected components of two link diagrams 
$D_1$ and $D_2$ such that 
$D_2$ is obtained from $D_1$ by smoothing at a crossing of $D_1$
is $0$ or $1$).
This contradicts the fact that $O_{+}(D)-1=n_{-}(D)$.
Therefore all negative crossings of $D$ are nugatory
and $D$ represents a positive knot.
\end{proof}

Corollary \ref{corollary:positivity} immediately follows from
Theorem \ref{theorem:positivity}.
\section{A new approach to estimate the Rasmussen invariant of a knot}
\label{section:A new approach to estimate the Rasmussen invariant of a knot}
Let $D$ be a diagram of a knot with $\Delta (D) \neq 0$.
Then Kawamura-Lobb's inequality may not be sharp as follow.
\begin{example}
Let $K$ be the pretzel knot of type $(3,-5,-7)$ 
and $D$ the standard pretzel diagram of $K$.
Then $\omega (D)=9, O(D)=14, O_{+}(D)=3$  and $O_{-}(D)=11$.
Therefore $\Delta (D)=1$ and $w(D)-O(D)+2O_{+}(D)-1=0$.
On the other hand, since $K$ is strongly quasipositive
\cite{Rudolph3},\footnote{In \cite{Rudolph3}, $K$ 
is denoted by $P(-3,5,7)$.
Our notation is the same as that in \cite{cromwell-book} and \cite{Kaw}.
} we obtain $s(K) =2g_{*}(K)=2g(K)=2$. 
Remark that $K$ is topologically slice but not smoothly slice.
\end{example}

We need a more shaper estimation to describe the Rasmussen invariant
of the pretzel knot of type $(3,-5,-7)$
in terms of its standard pretzel diagram.
Roughly speaking, there are two approaches
 to estimate or determine the Rasmussen invariant.
One of them is to compute the Khovanov homology by using a computer
and to use the spectral sequence which converges to Lee's homology.
The other is to use some formal properties of the Rasmussen invariant
(and the tau invariant). 
We propose a new and direct approach to estimate or
determine the Rasmussen invariant.
We briefly recall the definition of the Rasmussen invariant to explain this.
For a full explanation, see \cite{Rasmussen}.

Let $D$ be a diagram of a knot $K$ and $C_{Lee}^{*}(D)$ Lee's complex
(see \cite{Rasmussen} for the definition).
Then Lee \cite{Lee} proved that the homology group of $C_{Lee}^{*}(D)$ is 
independent of the choice of diagrams of $K$. 
Lee's homology of $K$, denoted by $H_{Lee}^{*}(K)$, is defined to be
the homology group of $C_{Lee}^{*}(D)$.
In addition, for a diagram $D$ of a knot $K$,
Lee \cite{Lee} associated two (co)cycles of $C_{Lee}^{*}(D)$,
denoted by $f_{o}$ and $f_{\bar{o}}$,\footnote{
In \cite{Rasmussen}, 
these are denoted by $s_{o}$ and $s_{\bar{o}}$ respectively} and proved that 
$[f_{o}]$ and $[f_{\bar{o}}]$ are a basis of $H_{Lee}^{*}(K)$, in particular,
that the dimension of  $H_{Lee}^{*}(K)$ is equal to two,
where [$\cdot$] denotes its homology class.
This basis is called canonical
 since the basis is determined 
up to multiple of $2^{c}$ for $K$ \cite{Rasmussen}, 
where $c$ is an integer.
 
Rasmussen \cite{Rasmussen} defined a filtration grading $q$ on a non-zero element of $C_{Lee}^{*}(D)$
(which induces a filtration on $C_{Lee}^{*}(D)$). 
Then a filtration grading $s$ on a non-zero element $[x]$ of $H_{Lee}^{*}(K)$ 
(which also induces a filtration on $H_{Lee}^{*}(K)$) 
is defined as follows. 
\[s([x])=\max\{q(y)|[x]=[y]\}. \]
Then the Rasmussen invariant of $K$, denoted by $s(K)$,  is defined to be $s([f_{o}])+1 (=s([f_{\bar{o}}])+1)$.

Since $s([f_{o}]) \ge q(f_{o})$ and $q(f_{o})=\omega(D) -O(D)$ (by the definition of $q$), we obtain $s(K) \ge \omega(D) -O(D)+1$.
This is the slice-Bennequin inequality for the Rasmussen invariant for $K$
 (see \cite{Plamenevskaya} and \cite{Shumakovitch}).
Theorem \ref{theorem:Kawamura-Lobb} implies that 
there exists a cycle $f$ such that $[f_{o}]=[f]$ and $q(f)=\omega(D) -O(D) +2O_{+}(D)-2$,
however, yet no one has succeeded to describe $f$ explicitly.
In \cite{Abe}, as a first step toward this, we describe 
a cycle $f_{1}$ with $[f_{o}]=[f_{1}]$ which gives the so-called 
shaper slice-Bennequin inequality for the Rasmussen invariant of a knot
\cite{Kawamura}
(which is stronger than the slice-Bennequin inequality and
weaker than the inequality of Kawamura and Lobb,
 see \cite{Kawamura2}).
In the future work, the graph $G(D)_{\Delta}$ is expected to
play an important role (see also \cite{Elliott}). 
We conclude this paper by giving the following problem.\\

\noindent
\textbf{Problem}. Let $D$ be the standard diagram of $P(3,-5,-7)$.
Find a cycle $f$ of $C_{Lee}^{*}(D)$ such that $[f_{o}]=[f]$ and
 $q(f)=\omega(D) -O(D) +2O_{+}(D)=1$.


\end{document}